\documentclass[12pt,a4paper]{amsart}
\theoremstyle{plain}

\usepackage[colorlinks]{hyperref}
\usepackage{enumerate, amssymb}

\advance\hoffset-5mm \advance\textwidth40mm
\input{diagrams.tex}
\diagramstyle[scriptlabels,height=8mm,width=8mm]

\def\bdi{\begin{diagram}}
\def\edi{\end{diagram}}

\theoremstyle{plain}

\newtheorem{thm}{Theorem}[section]
\newtheorem{cor}[thm]{Corollary}
\newtheorem{lem}[thm]{Lemma}
\newtheorem{prop}[thm]{Proposition}
\theoremstyle{definition}
\newtheorem{defi}[thm]{Definition}
\newtheorem{defis}[thm]{Definitions}
\newtheorem{conj}[thm]{Conjecture}
\newtheorem{conv}[thm]{Convention}
\newtheorem{nota}[thm]{Notation}
\newtheorem{rem}[thm]{Remark}
\newtheorem{rems}[thm]{Remarks}
\newtheorem{exa}[thm]{Example}
\newtheorem{exas}[thm]{Examples}
\newtheorem{prob}[thm]{Problem}
\newtheorem{probs}[thm]{Problems}
\newtheorem{ques}[thm]{Question}
\newtheorem{sit}[thm]{}

\newcommand{\Ker}{ \operatorname{{\rm Ker}}}

\renewcommand{\epsilon}{\varepsilon}

\def\and{\quad\mbox{and}\quad}

\newcommand{\C}{\ensuremath{\mathbb{C}}}

\newcommand{\N}{\ensuremath{\mathbb{N}}}

\newcommand{\bX}{{\bar X}}

\newcommand{\cO}{{\ensuremath{\mathcal{O}}}}

\renewcommand{\rho}{\varrho}

\def\bals#1\eals{\begin{align*}#1\end{align*}}
\def\bal#1\eal{\begin{align}#1\end{align}}

\def\PP{{\mathbb P}}

\renewcommand{\phi}{\varphi}

\newcommand{\bnum}{\begin{enumerate}}
\newcommand{\enum}{\end{enumerate}}
\renewcommand{\emptyset}{\varnothing}

\newcommand{\brem}{\begin{rem}}
\newcommand{\brems}{\begin{rems}}
\newcommand{\erem}{\end{rem}}
\newcommand{\erems}{\end{rems}}
\newcommand{\bprob}{\begin{prob}}
\newcommand{\eprob}{\end{prob}}
\newcommand{\bprobs}{\begin{probs}}
\newcommand{\eprobs}{\end{probs}}
\newcommand{\bques}{\begin{ques}}
\newcommand{\eques}{\end{ques}}
\newcommand{\bexa}{\begin{exa}}
\newcommand{\bexas}{\begin{exas}}
\newcommand{\eexa}{\end{exa}}
\newcommand{\eexas}{\end{exas}}
\newcommand{\bdefi}{\begin{defi}}
\newcommand{\edefi}{\end{defi}}
\newcommand{\bdefis}{\begin{defis}}
\newcommand{\edefis}{\end{defis}}
\newcommand{\bcor}{\begin{cor}}
\newcommand{\ecor}{\end{cor}}
\newcommand{\blem}{\begin{lem}}
\newcommand{\elem}{\end{lem}}
\newcommand{\bconv}{\begin{conv}}
\newcommand{\econv}{\end{conv}}
\newcommand{\bconj}{\begin{conj}}
\newcommand{\econj}{\end{conj}}
\newcommand{\bprop}{\begin{prop}}
\newcommand{\eprop}{\end{prop}}
\newcommand{\bthm}{\begin{thm}}
\newcommand{\ethm}{\end{thm}}
\newcommand{\bnota}{\begin{nota}}
\newcommand{\enota}{\end{nota}}
\newcommand{\bsit}{\begin{sit}}
\newcommand{\esit}{\end{sit}}
\newcommand{\be}{\begin{equation}}
\newcommand{\ee}{\end{equation}}
\newcommand{\bproof}{\begin{proof}}
\newcommand{\eproof}{\end{proof}}
\def\ba{\begin{array}}
\def\ea{\end{array}}

\begin{document}
\title[Analytic extensions of algebraic isomorphisms]{Analytic extensions of algebraic isomorphisms}

\author{S.\ Kaliman}
\address{Department of Mathematics,
University of Miami, Coral Gables, FL 33124, USA}
\email{kaliman@math.miami.edu}

\begin{abstract} Let $\Psi : X_1 \to X_2$ be an isomorphism of closed affine algebraic subvarities of $\C^n$ such that
$n > \max (2\dim X_1, \dim TX_1)$. We prove that $\Psi$ can be extended to a holomorphic automorphism of
$\C^n$. Furthermore, when $\Psi$ is an isomorphism of curves such an extension exists for every $n\geq 3$ even when
$\dim TX_1=n$.

\end{abstract}
\date{\today}
\maketitle

\thanks{
{\renewcommand{\thefootnote}{} \footnotetext{ 2010
\textit{Mathematics Subject Classification:}
14R20,\,32M17.\mbox{\hspace{11pt}}}}

\vfuzz=2pt
\thanks{}

%

\section*{Introduction}

Let $\Psi : X_1 \to X_2$ be an isomorphism of closed affine algebraic subvarities of $\C^n$. It was proven in \cite{Ka91}, \cite{Sr} that $\Psi$ is the restriction
of an algebraic automorphism of $\C^n$ provided $n> \max (2\dim X_1+1, \dim TX_1)$. The inequality
on $\dim TX_1$ cannot be improved \cite{Ka91} while the question whether $2 \dim X_1 +1$ is the optimal bound in this theorem remains open.
It is unknown, for instance, whether an isomorphism of two closed curves in $\C^3$  (even when these curves are isomorphic
to $\C$) can be extended to an algebraic automorphism of $\C^3$. However, any two algebraic\footnote{Precaution: the similar statement
does not hold for proper holomorphic embeddings of $\C$ into $\C^n$ for any $n \geq 2$, see \cite{FGR}, \cite{Ka92}.}
embeddings of $\C$ into $\C^3$ are equivalent up to a holomorphic coordinate substitution \cite{Ka92}. In this paper we get
the following improvements of the last statement.

\bthm\label{main} Let $\Psi : X_1 \to X_2$ be an isomorphism of closed affine algebraic subvarities of $\C^n$ and $n> \max (2\dim X_1, \dim TX_1)$.
Then $\Psi$ can be extended to a holomorphic automorphism $\tilde \Psi$ of $\C^n$.
\ethm 

\bthm\label{curve} Let $X_1$ and $X_2$ be closed affine algebraic curves in $\C^{n}, \, n \geq 3$, and
$\Psi : X_1 \to X_2$ be their isomorphism. Then $\Psi$ extends to a holomorphic automorphism of $\C^{n}$.

\ethm

In particular, in the case of an isomorphism of closed curves in $\C^3$ there are no analytic invariants that can prevent an extension to 
an algebraic automorphism. 

We are grateful to G. Freudenburg who shared with the author some of his unpublished results on algebraic extensions of
curve isomorphisms in $\C^3$ and whose questions stimulated this study.

\section{Semi-isomorphisms}

Our central technical tools will be the following notion of semi-isomorphism and also a weaker notion of pseudo-isomorphism
that will be used in Section 3.

\bdefi\label{semi} Let $\varphi : X \to Y$ be a holomorphic map of Stein spaces. We call $\varphi$ a semi-isomorphism
(resp. pseudo-isomorphism) if the following conditions hold

(a) $\varphi$ is a finite bimeromorphic map;

(b) for every $x \in X$ (resp. every smooth point $x \in X$) there exists a neighborhood $U \subset X$ of $x$ such that $\varphi |_U$ is an embedding;

(c) there exists a neighborhood $U_0$ of  the  singular  set $X_{\rm sing}$ of  $X$ such that $U_0=\varphi^{-1} (\varphi (U_0))$ and 
the restriction of $\varphi |_{U_0} : U_0 \to \varphi (U_0)$  is an embedding (resp. injection),
while the singular set $Y_{\rm sing}$ of $Y$ is a disjoint union of $\varphi (X_{\rm sing})$ and a discrete set of nodes (i.e. singular points $y\in Y$
such that the analytic germ of $Y$ at $y$ is a union of two smooth irreducible branches meeting transversally at $y$).

\edefi

We apply these notions mostly in algebraic case when $\varphi : X \to Y$ is a morphism of affine algebraic varieties over $\C$
(which will be the ground field throughout  the paper).
The main property of semi-isomorphisms we are going to exploit is given by the next straightforward fact.

\bprop\label{algebra} Let $\varphi : X \to Y$ be  a semi-isomorphism as in Definition \ref{semi}. 
Let $y_1, y_2 \ldots $ be the set of nodes of $Y$ not contained in  $\varphi (X_{\rm sing})$.
Suppose that $\varphi^{-1} (y_i) = \{ x_i',x_i'' \}$ for every $i$.
Then the algebra $\cO (Y)$ of holomorphic functions on $Y$
is naturally isomorphic to the subalgebra $R$ of the algebra $\cO (X)$ of holomorphic
functions on $X$ such that $f \in R$ if and only if $f(x_i')=f(x_i'')$ for every $i=1, \ldots , k$.  

\eprop

To describe semi-isomorphisms geometrically in the algebraic  situation we need the following.

\bnota\label{chord}
For  $m\leq n$  the space of linear surjective maps (projections) from $\C^n$ to $\C^m$ will be denoted by $L_{n,m}$.
This notion depends on the choice of coordinate
system in $\C^n$, i.e. on the choice of an embedding of $\C^n$ into $\PP^n =\C^n\cup H$ where $H$ is the hyperplane
at infinity. 
For every affine subvariety $X$ of $\C^n$ denote by $\bX$ its closure in $\PP^n$ and by $X_H$ the variety
$\bX \cap H$.  
For every $x \in X$ one can treat the Zariski tangent space $T_xX$ as a linear subspace of $\C^n$.
With this interpretation we consider the variety $T'X =\{ x+v | x \in X, v \in T_xX \}$. We let
$Ch (X)$ be the variety of chords of $X$  (recall that a chord is a line that passes at least through two points of $X$)
and  $CX$ be the closure (in $\C^n$) of $\bigcup_{\ell \in Ch (X)} \ell$.
Note that there is a natural morphism $\pi : Ch (X) \to (CX)_H$ such that for general $z \in (CX)_H$ the preimage $\pi^{-1} (z)$
consists of parallel chords of $X$ whose direction is determined by $z$.

\enota

We need also the following.

\blem\label{condition} Let $X$ be a closed proper affine algebraic subvariety of $\C^n$such that  $n\geq 2\dim X +1$. 
Then the coordinate system in $\C^n$ can be chosen so that for any pair of general points $x_1,x_2$ in $X$ (not necessarily in the same irreducible component
\footnote{Recall that a point of 
an algebraic variety $X$ is called general if it is not contained in a proper closed subvariety $Z$ of $X$. We require additionally that such $Z$
does not contain any irreducible component of $X$.})
the following condition holds

\noindent ${\rm (\# )} \, \, \, $ \hspace{1cm} 
the vector space generated by  $T_{x_1}X, T_{x_2}X$ and the vector $\overrightarrow{x_1,x_2}$ has dimension
$\dim T_{x_1}X + \dim T_{x_2}X+ 1$ (in particular, it is  $2 \dim X +1$ when $x_1$ and $x_2$ belong to components of $X$ with the largest
dimension).

\elem

\bproof 

Let $A$ be a finite set of general points of $X$ such that every irreducible component of $X$ contains at least
two points from $A$.  By \cite[Theorem 4.14 and Remark 4.16]{AFKKZ} there exists an algebraic automorphism $\beta$ of $\C^n$
that fixes every point $a\in A$ and has
prescribed values of $(\beta)_*|_{T_a\C^n}$ in $SL(T_{a}\C^{n})$ for each $a \in A$. Choosing appropriate 
collection $\{ (\beta)_*|_{T_a\C^n} | a \in A \}$ and replacing $X$ by $\beta (X)$
we get Condition ($\#$).
\eproof

\brem\label{rem1.1}  (1)  Let $ \varphi\in L_{n,m}$ be such that  $\overrightarrow{x_1,x_2} \in \Ker \varphi$ for $x_1, x_2$ as in Lemma \ref{condition}. Then
Condition ($\#$) implies that
$\varphi (T_{x_1}X) \cap \varphi (T_{x_2}X) = \{ 0 \}$.   
 
 (2) One can apply \cite[Theorem 4.14 and Remark 4.16]{AFKKZ} in the proof of this Lemma \ref{condition} because $\C^n$ is 
 so-called flexible variety. Let $Y$ be a closed subvariety of $\C^n$ with $\dim Y \leq n-2$. Then $\C^n \setminus Y$ is still flexible (e.g., see \cite{FKZ13}).
 In particular, if $Y$ does not contain any irreducible component
 of $X$ then the automorphism of $\C^n$ that induces a new coordinate system in Lemma \ref{condition} can be chosen so that 
 its restriction to $Y$ (and even to any given infinitesimal neighborhood of $Y$) is identical \cite[Theorem 1.6]{FKZ13}.
 
 (3) The similar proof shows that Lemma \ref{condition} is valid in the case when $X$ is a closed proper analytic subset of $\C^n$ with a finite number of
 irreducible components. A bit more complicated argument implies that it remains true even for an infinite number of components (but we shall not
 need this fact later). 
\erem

\bconv\label{conv1.1} Further in this section we suppose that $X$ is a closed proper affine algebraic subvariety of $\C^n$ 
such that $n\geq 2\dim X+1$ and $\dim TX \leq 2 \dim X$ where $TX$ is the Zariski
tangent bundle of $X$. We continue to treat every tangent space $T_xX$ as a natural subspace of $\C^n$ and suppose that
Condition ($\# $) of Lemma \ref{condition} holds. 
\econv

\blem\label{lem1.1} Let Convention \ref{conv1.1} hold.
Then $\dim CX = 2\dim X +1$ and therefore $\dim (CX)_H = 2\dim X$. Furthermore,
there is a subvariety $ (CX)_H' \subset  (CX)_H$ of dimension at most $2\dim X -1$ such that for every 
$z \in  (CX)_H\setminus  (CX)_H'$

{\rm (i)} the preimage $\pi^{-1}(z)$ is at most finite;

{\rm (ii)} for each $\ell \in \pi^{-1} (z)$ every pair of distinct points $x_1$ and $x_2$ in $\ell \cap X$
satisfies Condition ($\#$) \footnote{That is, for general  $\varphi\in L_{n,n-1}$
and every pair of points $x_1, x_2 \in X$ such that  $\overrightarrow{x_1,x_2} \in \Ker \varphi$ one has  Condition ($\#$).};

{\rm (iii)}  $\ell$ does not contain points from $X_{\rm sing}$ and $\ell$ is not tangent to $X$ at any smooth point;

{\rm (iv)}  $\ell$ meets $X$ exactly at two points.

\elem

\bproof Consider general points $x_1$ and $x_2$ in an irreducible component $D$ of $X$ such that $\dim D = \dim X$.
Let $U_i$ be a small Euclidean neighborhood of $x_i$ in $A$, i.e. $U_i$ corresponds to a neighborhood of the origin
in $T_{x_i}D$. Condition ($\#$) implies that the union of chord through points of $U_1$ and $U_2$
is of dimension $2 \dim X+1$ which yields the first statement.

The union $B$ of chords containing points from $X_{\rm sing}$
has dimension at most $ \dim X_{\rm sing} + \dim X +1 \leq 2\dim X$. That is, $\dim B_H < \dim H$. Similarly by Convention
\ref{conv1.1} $\dim (T'X)_H < 2\dim X$.  Choosing $(CX)_H'$ so that  $B_H$ and $(T'X)_H$ are contained in $(CX)_H'$ we get (iii). 
In particular every $\ell$ is a line that meets $X$ transversally at the intersection points.

Suppose that $\pi^{-1} (z)$ is not finite. That is, assuming that  $x_1, x_2 \in \ell$ and $U_1, U_2$ are as before we can 
find biholomorphic analytic subsets $W_i\subset U_i$ of positive dimensions such that every $w_1 \in W_1$ and its image $w_2$ in $W_2$
are joined by a chord from $\pi^{-1} (z)$. This implies that for $\varphi$ from Remark \ref{rem1.1} (1)
$\varphi (T_{x_1}X) \cap \varphi (T_{x_2}X)$ contains $\varphi (T_{x_1} W_1) \simeq \varphi (T_{x_2} W_2)$. However
the points $x_1$ and $x_2$ are general since $z$ is general. This contradicts
Condition ($\#$) from Convention \ref{conv1.1}. That is, $W_i$ is zero-dimensional and we have (i).
By the same reason we have (ii).

Note that the number of points in $\ell \cap X$ is constant for general $z \in (CX)_H$ and these points depend
continuously on $z$. (Indeed, there is a natural morphism from the complement to the diagonal in $X \times X$ into $(CX)_H$
and the statement is a consequence of the semi-contnuity theorem \cite{Ha}.) Assume that
$\ell$ contains at least three points $x_1, x_2$ and $x_3$ in $X$, i.e. the vectors  $\overrightarrow{x_1x_2}$ and $\overrightarrow{x_1x_3}$ are parallel.
Let $U_i$ be a neighborhood $x_i$ as before. Choose a sequence of points $x_{2i}$ in $U_2$ convergent to $x_2$ 
and consider the chords through $x_1$ and $x_{2i}$. By semi-continuity and continuous dependence on $z$ there must
be a point $x_{3i} \in U_3$ on the same chord as $x_1$ and $x_{2i}$. However this implies that for $\psi \in L_{n, n-1}$ with $\ell \in \Ker \psi$ 
one would have again  $\psi (T_{x_2}X) \cap \psi (T_{x_3}X) \ne \{ 0 \}$
contrary to Condition ($\#$). Hence we have  (iv).

\eproof

The following is the main result of this section.

\bthm\label{criterion}
Let $X$ be a closed affine subvariety of $\C^n$ with Convention \ref{conv1.1} valid, and let
$2 \dim X \leq m < n$. Suppose that $\varphi$ is a general element of  $L_{n,m}$.
Then $\varphi |_X : X \to Y:=\varphi (X) \subset \C^m$ is a semi-isomorphism.
\ethm

\bproof  
Since for $m=2 \dim X +1$ and general $\varphi \in L_{n,m}$ the map
$\varphi |_X : X \to Y$ is an isomorphism onto a closed subvariety of $\C^m$ (e.g., see \cite{Ka91})
we can suppose that $n=2 \dim X +1$ and $m=2 \dim X$.
Then by Lemma \ref{lem1.1} (i) $\varphi |_X : X \to Y$ is bijective in the complement to a finite set for general $\varphi$.
That is,  such a $\varphi$ is birational.

Since $V= \Ker \varphi$ is now a line,
$V_H$ is  a general point in $H$.   Then  $V_H \cap X_H=\emptyset$ because $\dim X_H = \dim X -1$.
Hence  $\varphi |_X : X \to Y$ is a proper morphism, and $Y$ is closed in $\C^m$. Furthermore
being proper and quasi-finite, $\varphi |_X$ is finite by Grothendieck's theorem.

Similarly $V_H \cap (T'X)_H=\emptyset$ for a general $\varphi$ since $\dim TX \leq 2\dim X$ and thus $\dim (T'X)_H \leq 2 \dim X -1$.
Hence for every $x \in X$ the linear map $\varphi_* : T_xX' \to T_{\varphi (x)}Y'$ is an isomorphism  where $X'$ is the analytic
germ of $X$ at $x$ and $Y'=\varphi (X')$.  This implies that for some neighborhood $U \subset X$
the restriction $\varphi |_U$ is an \'etale map (e.g, see \cite[Proposition 7]{Ka91}). 

Thus we have Conditions 
(a) and (b) from Definition \ref{semi} while Condition (c)  follows from Lemma \ref{lem1.1} (ii), (iii), and (iv).

\eproof

\brem\label{rem1.2} 
(1)  Since $\varphi \in L_{n,m}$ in Theorem \ref{criterion} is general we can suppose that
any forgetting projection $\C^n \to \C^m$ to a coordinate subspace of dimension $m$ can serve as $\varphi$ in this Theorem. 

(2) Since any semi-isomorphism
does not distinguish only a finite number of points and any general linear function separates elements of finite sets we
can also suppose that any forgetting projection $p: \C^n \to \C^{m+1}$ yields an isomorphism between $X$ and
$p(X)$. 

\erem

Furthermore, one can see that the proof yields a stronger statement.

\bprop\label{prop1.8} Let $X$ and $Y$ be as in Theorem \ref{criterion}.
If  $l= 2\dim X$ and $M$ is a subvariety of $L_{n,l}$ of dimension
at least $l+1$ such that for a general $\varphi \in M$ and  any points $x_1, x_2 \in X$ with $\varphi (x_1)=\varphi (x_2)$
Condition ($\#$) from Lemma \ref{condition} holds, then the restriction $\varphi |_X : X \to Y$ is a semi-isomorphism.
\eprop

\bproof
Indeed,  for $1 \leq m \leq n$ and every $\psi  \in L_{n,m}$ let us treat $\psi (\C^n )$ as orthogonal complement to $\Ker \psi$  (i.e. as a subspace  of $\C^n$).
Then one can fix a general $\psi \in L_{n, l+1}$ such that  $\psi |_X : X \to \psi (X)$ is an isomorphism onto the closed subvariety
$\psi (X)$ of $\C^{l+1}$ and for general $\varphi \in M$  the restriction of $\psi$ to $\varphi (\C^n ) \subset \C^n$ is injective.
We can present now $\psi \circ \varphi : \C^n \to \C^l$ as a composition $\varphi' \circ \psi : \C^n \to \C^l$
where $\varphi' \in L_{l+1,l}$, and replace $X$ by $\psi (X)$ and $M$ be the family $M'$ that consists of such $\varphi'$. 
Generality of $\psi$ implies that its restriction to the subspace generated by  $T_{x_1}X, T_{x_2}X$ and the vector $\overrightarrow{x_1,x_2}$
is injective and therefore
the modified Condition ($\#$) does not suffer under this change from $M$ to $M'$  which yields Condition (c) from Definition \ref{semi}
while Conditions (a) and (b) follow as before from the fact that $\dim M' =l+1> \dim (CX)_H'$.
\eproof

To describe a specific submanifold $M$ needed later we introduce the following.

\bnota\label{nota1.2} Let $\Psi : X_1 \to X_2$ be an isomorphism of  closed affine algebraic subvarieties of $\C^n$ where $n=2 \dim X_1 +1$.
Suppose that Convention \ref{conv1.1} is true for the embeddings $X_i \subset \C^n, \, \, i=1,2$, and the restriction of forgetting projections satisfy Remark \ref{rem1.2} (1).
Consider the graph $\Gamma$ of $\Psi$ in $\C^{2n}$ with natural coordinates $(\bar z, \bar w)=(z_1, \ldots , z_n, w_1, \ldots , w_n )$
and the forgetting projections $\rho_{m,k} : \C^{2n} \to \C^{m+k}$ given by
$$ (\bar z, \bar w) \to (z_1, \ldots , z_m, w_1, \ldots , w_k).$$ 
\enota

\bprop\label{prop1.3} Let  Notation \ref{nota1.2} hold and $l= 2\dim X_1$. Suppose that $M$ is the subvariety of $L_{2n, l}$
that consists of all projections of form $\varphi = (\varphi_1, \varphi_2)$ where $\varphi_1$ depends on $\bar z$ only while
$\varphi_2$ depends on $\bar w$ only.Then there is a polynomial coordinate substitution in $\C^{2n}=\C^n \oplus \C^n$
of form $\alpha = (\alpha_1, \alpha_2)$ where $\alpha_i$ is an algebraic automorphism of $\C^n$ such that after this substitution
the restriction $\varphi |_\Gamma : \Gamma \to \varphi (\Gamma)$ is a semi-isomorphism onto a closed subvariety of $\C^l$ for any general $\varphi \in M$.
In particular, one can assume (by Remark \ref{rem1.2}) that the restriction to $\Gamma$ of each $\rho_{2n,l}$ is a semi-isomorphism and 
of each $\rho_{2n,n}$ is an isomorphism.
\eprop

\bproof  Consider general points $a_1$ and  $a_2$ in $\Gamma$. Suppose that $\varphi \in M$,
i.e.  $\varphi = (\varphi_1, \varphi_2) : \C^{2n}=\C^n\oplus \C^n \to \C^l = \C^m \oplus \C^{l-m}$. We can assume $m \geq \dim X$ (if not,
switch $m$ and $l-m$). Let  $a_{ij}=\varphi_i(a_j)$ (in particular, $a_{i1}$ and $ a_{i2}$
are general points in $X_i$) and $V_{ij}=(\varphi_i)_* (T_{a_j} \Gamma)$.  Generality and Remark \ref{rem1.2} imply that $\dim V_{1j}=\dim T_{a_{1j}}X_1$ for every $j$.
Choose every $\alpha_i$ so that it preserves $a_{i1}$ and $a_{i2}$.
Since $n \geq 2$, by \cite[Theorem 4.14 and Remark 4.16]{AFKKZ} $\alpha_i$ can be also chosen
with prescribed values of $(\alpha_i)_*$ in $SL(T_{a_{ij}}\C^{n})$ for each $j$. In particular, we can make this choice so 
that the space generated by $ \alpha_{1*} (V_{11})$,  $ \alpha_{1*} (V_{12})$, and the vector $\overrightarrow{a_{11}a_{12}}$
has dimension $\dim T_{a_{11}}X_1+ \dim T_{a_{12}}X_1 +1$.  That is, Condition ($\#$) from Lemma \ref{condition} holds for points $a_1$ and $a_2$.
Since $\varphi_1$ (resp. $\varphi_2$) is a general element of $L_{n,m}$
(resp. $L_{n,l-m}$) and $a_{ij}$ is a general point of $X_i$ we can suppose that after a perturbation $\varphi_i (a_{i1})=\varphi_i(a_{i2})$  for $i=1,2$,
i.e.  $\varphi (a_{1})=\varphi (a_{2})$. 
Thus  we are done now by Remark \ref{rem1.2}.


\eproof

\section{The proof of Theorem \ref{main}}

\bprop\label{prop2.1}  Let $\Psi : X_1 \to X_2$ be an isomorphism of closed affine algebraic subvarieties of $\C^n$ where $n =2 \dim X_1+1$.
Suppose that $\varphi \in L_{n,n-1}$ is such that $\varphi (X_1)=\varphi (X_2)=:Y$ is closed in $\C^n$ and $\varphi |_{X_i} : X_i \to Y$ is a semi-isomorphism
for $i=1,2$. Then there exists an extension of $\Psi$ to a holomorphic automorphism $\tilde \Psi : \C^n \to \C^n$.

\eprop

\bproof Let $\bar z = (z_1, \ldots , z_n)$ be a coordinate system on $\C^n$. Without loss of generality one can suppose that
the coordinate form of $\varphi$ is $\bar z \to (z_1, \ldots , z_{n-1})=:\bar \zeta$. Denote by $f_i$ the restriction of the function $z_n$ to $X_i$.
Then $f_i$ can be treated as a rational function on $Y$ with possible indeterminacy points in $Y_{\rm sing} \setminus S$
where $S= \varphi ((X_1)_{\rm sing})= \varphi ((X_2)_{\rm sing})$ (the last equality holds since otherwise $\Psi$ is not an isomorphism).
Consider a holomorphic automorphism $\alpha$ of $\C^n$ given by $\bar z \to (\bar \zeta , e^{g(\bar \zeta )}z_n)$ where $g$ is a polynomial in $\bar \zeta$. 
Our aim to choose $g$ so that $h:=f_1-e^{g(\bar \zeta)}f_2$ is a holomorphic function on $Y$.

Let $y \in Y_{\rm sing} \setminus S$ and $\varphi^{-1} (y) \cap X_i = \{ x_i', x_i'' \}$. Suppose that $a_i' =f_i (x_i')$ and $a_i''=f_i(x_i'')$.
Note that $a_i' \ne a_i''$ since otherwise $x_i'=x_i''$ is a singular point of $X_i$ contrary to the assumption.
Choose $g$ so that $$e^{g(y)}=\frac{a_1'-a_1''}{a_2'-a_2''}$$ for every $y \in Y_{\rm sing} \setminus S$. Then $h(x_1')=h(x_2'')$
and therefore $h$ can by viewed as a holomorphic function on $Y$ by Proposition \ref{algebra}. We denote by the
same symbol a holomorphic extension of $h$
to $\C^{n-1}$.
Let $\beta$ be the holomorphic automorphism of $\C^n$ given by $\bar z \to (\bar \zeta, z_n + h(\bar \zeta))$.
Then the composition $\alpha^{-1} \circ \beta$ is the desired extension of $\Psi$.
\eproof

\brem\label{rem2.5}
(1) The same argument shows that an extension to a holomorphic automorphism in Proposition \ref{prop2.1} exists 
in the case when $\Psi : X_1 \to X_2$ is a biholomorphic map of closed analytic subsets of $\C^n$ and $Y$ is a closed analytic subset
of $\C^{n-1}$.

(2) If $\varphi$ is only a pseudo-isomorphism then a priori the meromorphic function $h=f_1 -e^{g(\zeta )}f_2= f_1-f_2 + (1-e^g(\zeta ) f_2$ constructed in the proof
is holomorphic only on $Y\setminus \varphi (X_{\rm sing})$. However, suppose that $X_{\rm sing}$ is, say, discrete and for every $x \in X_{\rm sing}$ there exists a natural
$k=k (x)$ such that every germ of a continuous meromorphic function, that vanishes at $x$ with multiplicity at least $k$, is in fact holomorphic at $x$.
If it is also known that $f_1-f_2$ is holomorphic at points of $\varphi (X_{\rm sing})\subset Y$
then requiring additionally that $g$ vanishes on $\varphi (X_{\rm sing})$ with sufficiently large multiplicity one gets
$h$ holomorphic on the entire $Y$. In particular in this case an extension of $\Psi$ to a holomorphic automorphism again exists.
\erem

\bdefi\label{defi2.1} Let Notation \ref{nota1.2} hold but we allow $n \geq  2 \dim X_1 +1$.
We say that the triple $(\Psi , X_1, X_2 )$ is admissible with respect to the coordinate system $(\bar z , \bar w )$ if

(i) $n >l:= \max (2\dim X_1, \dim TX_1)$;

(ii)  for everhy $m=0, 1, \ldots , l$ the set $Z_m= \rho_{m, l-m}(\Gamma )$ is a closed affine algebraic subvariety of $\C^l$;

(iii) the restriction  $\rho_{m, l-m}|_{\Gamma} : \Gamma \to Z_m$ is a semi-isomorphism for $m=0, \ldots , l$.

(iv) the restriction  $\rho_{m, l+1-m}|_{\Gamma} : \Gamma \to W_m:=\rho_{m, l+1-m}(\Gamma )$ is an isomorphism for $m=0, \ldots , l+1$.

\edefi

\bprop\label{prop2.2} Let a triple $(\Psi , X_1, X_2 )$ be admissible. Then $\Psi$ can be extended to a holomorphic
automorphism $\tilde \Psi : \C^n \to \C^n$.

\eprop

\bproof If $\dim TX_1 \geq 2\dim X_1+1$ then $\tilde \Psi$ can be chosen algebraic by \cite{Ka91}. Thus we suppose that 
$\dim TX_1 \leq 2\dim X_1$. By the same reason we can suppose that $n=2\dim X_1 +1=l+1$.
Consider $\C^n$ equipped with a coordinate system $(u_1, \ldots , u_n)$.
Let $W_m'$ be the image of $\Gamma$ under the morphism $\C^{2n} \to \C^n$ given by $$(\bar z , \bar w) \to (z_1, \ldots , z_{m-1}, w_1, \ldots , w_{n-m}, z_{m}).$$
Note that all $W_m$ are closed subvarieties of $\C^n$ by Condition (ii) of Definition \ref{defi2.1}. The same is true for $W_m'$ since
$W_m$ and $W_m'$ are isomorphic under the automorphism $\beta_m : \C^n \to \C^n$ that switches coordinates $u_m$ and $u_n$. 
Note also that $W_0=X_2$ while $W_n=X_1$.
By Condition (iv) of Definition \ref{defi2.1} there are isomorphisms $\psi_m: W_{m-1} \to W_m$ such that $\Psi = \psi_1 \circ \cdots \circ \psi_{n}$.
Consider the projection $\varphi : \C^n \to \C^l$ given by $(u_1, \ldots , u_n) \to (u_1, \ldots , u_l)$. Note that $\varphi (W_{m-1})=\varphi (W_m')$
and by Condition (iii) of Definition \ref{defi2.1} the restrictions of $\varphi$ to $W_m$ or to $W_m'$ yield semi-isomorphisms.
By Proposition \ref{prop2.1} there is a holomorphic automorphism $\alpha_m : \C^n \to \C^n$ that extends the isomorphism $\beta_m \circ \psi_{m}$.
Hence $\beta_m\circ \alpha_m: \C^n \to \C^n$ is the extension of $\psi_m$ and  $(\beta_1\circ \alpha_1) \circ \cdots \circ (\beta_n\circ \alpha_n)$ is the 
desired extension of $\Psi$.
\eproof  

\subsection{Proof of Theorem \ref{main}.} By \cite[Theorem 1]{Ka91} we can suppose that $\dim TX_1 \leq 2 \dim X_1$.
Let $l=2\dim X_1 +1$. Since for a general $\varphi \in L_{n,l}$ the restrictions $\varphi |_{X_i} : X_i \to \varphi (X_i), \, \, i=1,2$ are isomorphisms 
onto closed subvarieties of $\C^l$ we can suppose that $n=2\dim X_1 +1$. Applying algebraic automorphisms $\alpha_i$ of $\C^n$ to $X_i, \, \.  i=1,2$
we can suppose that the last statement of Proposition \ref{prop1.3} holds. In particular, the triple $(\Psi , X_1, X_2)$ is admissible.
Hence by Proposition \ref{prop2.2} $\Psi$ can be extended to a holomorphic automorphism $\tilde \Psi$ of $\C^n$ which is the
desired conclusion. \hspace{2.7cm} \hspace{3.4cm} $\square$

\section{The case of curves}

It is interesting to find out whether the restriction on $\dim TX_1$ in Theorem \ref{main} can be made weaker; say, whether in the
case of isolated singularities 
the inequality $n \geq 2 \dim X_1 +1$ yields the existence of an extension to a holomorphic automorphism. At present we do not
know a complete answer except for the one-dimensional case.
Recall that  there are isomorphic
curves in $\C^3$ whose isomorphism does not extend to a algebraic automorphism\footnote{An example of such curves is given
by the images of the morphisms $\rho_i : \C \to \C^3, \,$ where 
$\rho_1(t)=(t^7,t^{11},t^{13}) \, \, \, {\rm and} \, \, \, \rho_2(t)=(t^{7+14},t^{11},t^{13})$ with an isomorphism between them given
by $(x,y,z) \to (x+x^2,y,z)$. Actually, the argument
in \cite{Ka91} shows that there is no extension of this isomorphism to a holomorphic automorphism with constant Jacobian.} of $\C^3$ \cite[Example]{Ka91} but we shall
show that an extension to a holomorphic automorphism always exists.

\bconv\label{conv3.1} We study below affine algebraic subvarieties of $\C^{n}$ equipped with coordinates
$(x,y, \bar z)$ where $\bar z = (z_1, \ldots , z_{n-2})$ . For the sake of notation but without loss of generality we suppose further that $n=3$ and $(x,y,z)$ is a coordinate
system on $\C^3=\C^{n}$.

\econv

We shall need later holomorphic automorphisms of $\C^{n}$ preserving some coordinates and having prescribed jets at a finite
collection of points\footnote{Holomorphic automorphisms with prescribed jets at a finite set were constructed
in \cite{For} but the condition on preservation of coordinates does not allow a direct use of that result.}.

\bprop\label{prop3.1} Let $S =\{ s_1, \ldots , s_k \}$ be a finite subset in $\C^{3}$ with coordinates
$(x,y,z)$ and 
$s_i=(x_i,y_i, z_i)$.  Suppose that $ z_i  \ne z_j$ for $i \ne j$. Then for every natural number $m\geq 1$
and a collection of polynomials $f_i(x, y, z), \, i=1, \ldots , k$ with $f_i(0,0, 0 )=0$ there exists a holomorphic automorphism
$\alpha$ of $\C^{3}$ identical on $S$ and such that for every $i=1, \ldots , k$
the Taylor series of $\alpha$ at $s_i$ is of form $\alpha (x,y, z)=$
\[ (x,y, z) + (c_i\Delta x_i + f_i(\Delta x_i, \Delta y_i, \Delta z_i) + M^{m}(\Delta x_i,\Delta y_i, \Delta z_i),  M^{m}(\Delta x_i,\Delta y_i, \Delta z_i), \bar 0) \leqno{(1)} \]\label{for3.2}
where $M$ is the maximal ideal of $\C [x,y, z]$ corresponding to the origin, $c_i \in \C^*$, $\Delta x_i =x-x_i, \Delta y_i = y-y_i$, and $\Delta z_i= z -z_i$.

\eprop

\bproof Consider  a polynomial $h(y, z)$ such that $h(y_i, z_i) = \ln c_i$ and $h(\Delta y_i, \Delta z_i) -\ln c_i \in M^{m}(0,\Delta y_i, \Delta z_i)$ for every $i=1, \ldots , k$.
Then using composition with holomorphic automorphism $(x,y, z) \to (e^{h(y, z )}y, z)$ one can always suppose that each $c_i=1$.
Furthermore, applying triangular automorphism of form $(x,y,z) \to (x+p(y,z),y,z)$ we can suppose that each $f_i(x,y,z)$ is divisible by $x$.
Another simplifying observation is that it suffices to prove the statement in the case of $f_i(x,y,z)=x^lg(y,z)$ where $1\leq l \leq m$.
Indeed, if this is true then taking composition of such  holomorphic automorphisms with $l$ changing from 1 to $m$
one gets Formula (1) in general case.
 
Let $q_{l,m}$  be the sum of the first $m+1$ terms of the Taylor series for $1/(1+lx^{l-1})$ at the origin, i.e. the
Jacobian on the map $(x,y) \to (x+x^l, y(1+q_{l,m}(x))$ is $1+M^m(x,y, 0)$. By \cite[Theorem 4.14]{AFKKZ} there exists a polynomial automorphism 
$\beta_0=(\beta_0', \beta_0''): \C^2 \to \C^2$ whose $m$-jet at the origin coincides with $(x,y) \to (x+x^l, y(1+q_{l,m}(x))$.  Let $r_1( z)$ and $r_2 ( z)$ be polynomials
such that $r_1(z_i)=x_i$ and $r_2( z_i)=y_i$ for each $i$. 
Consider the automorphism $\beta$ of $\C^{3}$ given by $$\beta (x,y, z )= (r_1( z)+ \beta_0'(x-r_1( z)), r_2( z)+ \beta_0'' (y-r_2( z)), z).$$
Let $q(x)$ be the sum of the first $m$ terms of the Taylor series of $\ln q_{l,m}$ at the origin. Note that the composition $\beta_1$ of $\beta$ with
the holomorphic automorphism $(x,y, z) \to (x, e^{-q(x-r_1( z))}y, z)$ satisfies Formula (1) with $f_i(x,y, z)=x^l$ for each $i$.
Let $g(y, z)$ be a polynomial such that $g(y_i , z_i)=0$ for every $i$
and the Taylor series of $g$ at $(y_i, z_i)$ is of form $g_{i}(y-y_i, z - z_i)+ M^m(0,y-y_i, z - z_i)$
where $\{ g_{i}| \, i=1, \ldots ,k \}$ is a given collection of polynomials. Then the composition $\gamma_l$ of
$\beta_1$ with the holomorphic automorphism $(x,y, z) \to (e^{g(y,z)}x,y ,  z)$ satisfies Formula (1) with $f_i(x,y,z)=x^lg_{i}(y,z), \, i=1, \ldots , k$ and we are done.
\eproof

\blem\label{lem3.1}
Let $X$ and $Y$ be the germs of algebraic curves at the origin $o$ of $\C^{3}$ and $M$ be as in Proposition \ref{prop3.1}.
Suppose that $\psi : X \to Y$ is a bijective morphism. Then there exists a natural $m$ such that $M^m|_X \subset \psi^* (\cO (Y))$.
\elem

\bproof Suppose first that $Y$ is irreducible, i.e. it is the image of the germ $(\C , 0)$ under an injective holomorphic
map $\lambda : (\C, 0) \to (Y, o)$. 
Let $R =\lambda^* (\cO (Y))$
and let $S$ be the similar ring for $X$, i.e. $S$ can be viewed as a finitely generated $R$-module
containing $R$ as a submodule.
The set of zero multiplicities of functions from $R$ at $0 \in \C$ form a semi-group  $G$ 
and, furthermore,  there exists a natural $m$ such that $G$ contains $\{ n \in \N | n\geq m \}$\footnote{Since the author does not
know a reference to this well-known fact we sketch a proof. Let $k \in \N$ be the greatest common divisor of elements of $G$, i.e.
for some $m>0$ this semi-group contains $\{ nk | n \in \N , n \geq m/k  \}$. Assume that $k \geq 2$. There must be a function
$f \in \cO (Y)$ such that Taylor series $f \circ \lambda (t) =\sum_{i=1}^\infty a_it^i$ contains a nonzero term $a_jt^j$ with $i$ non-divisible by $k$ since otherwise the map
$\lambda$ is not injective. Replacing $f$ with its power one can suppose that $a_1=\ldots =a_m=0$ and still have some $a_i \ne 0$
with $i$ non-divisible by $k$. Hence adding to $f$ an element of $R$ we can get a function from $R$ whose zero multiplicity is not divisible by $k$. A contradiction. }.
 
 Suppose that $I$ is the ideal of $S$ induced by   $M^m|_X$. By the Krull theorem (e.g., see \cite{Ei}) $\bigcap_{k=1}^\infty I^k = 0$.
 Hence   $\bigcap_{k=1}^\infty I^k /(R\cap I)= 0$.
Then the above property of the semi-group $G$ implies that for every $k>0$ and $f \in I$ there exists $g\in R \cap I$ such that
$f+g \in I^k$. That is, the image of $f$ in $I/(R\cap I)$ is zero which 
implies the desired conclusion in the irreducible case.

Suppose now that $Y$ is not irreducible. Say, it consists of two components  $Y'$ and $Y''$. Let $m'$ and $m''$ play the same
role for this components as $m$ for $Y$ in the irreducible case. Let $g'$ be a holomorphic function that vanishes on $Y''$ but not 
on $Y'$ while $g''$ plays the opposite role. Let $n'$ be the zero multiplicity of $g'$ on $Y'$ and $n''$ the zero multiplicty
of $g''$ on $Y''$. The irreducible case implies now that any function on $Y$ that vanishes on $Y''$ (resp. $Y'$) and
has zero multiplicity  $n'+m'$ on $Y'$ (resp. $n''+m''$ on $Y''$)  is contained in $R$. Taking the sum of such functions
we get the desired conclusion with $m=\max (n'+m', n''+m'')$.

\eproof

\bprop\label{prop3.2} Let $\varphi : \C^3 \to \C^2_{y,z}$ be the forgetting projection, $\Psi : X_1 \to X_2$
be a isomorphism of closed affine algebraic curves in $\C^3$ such that $\varphi |_{X_1} = \varphi \circ \Psi$
and $\varphi |_{X_i} : X_i \to Y:=\varphi (X_i)$ is a pseudo-isomorphism. Suppose also that
we have $(X_1)_{\rm sing}=(X_2)_{\rm sing}=:S$. Then there exists an
extension of $\Psi$ to a holomorphic automorphism of $\C^3$.
\eprop

\bproof By assumption the coordinate form of $\Psi$ is $(x,y,z)\to (F(x,y,z),y,z)$. 
Making a linear substitution in $\C^2_{y,z}$ one can suppose that the restriction of
the forgetting projection $\C^3 \to \C_z$ to $S:=(X_i)_{\rm sing}$ is injective and since $\Psi$ is an isomorphism
at every $s_i=(x_i,y_i,z_i) \in S$ the Taylor series of $F$ is of form $F(x,y,z)=c_ix_i +f_i(\Delta x_i,\Delta y_i,\Delta z_i) +M^m(\Delta x_i,\Delta y_i,\Delta z_i)$
where $m$ is a sufficiently large natural number, $c_i, f_i, \Delta x_i,\Delta y_i,\Delta z_i$, and $M$ are as in Proposition \ref{prop3.1}.
Hence applying a holomorphic automorphism as in Proposition \ref{prop3.1} to $X_1$ we can suppose that $c_i=1$ and  $f_i\equiv 0$ for every $s_i \in S$.  
Consider $h = x|_{X_1} - x|_{X_2}$ as a meromorphic function on $Y$. Note that by Lemma \ref{lem3.1} $h$ is holomorphic at points $\varphi (S)$.
Now the desired conclusion follows from Remark \ref{rem2.5} (2).
\eproof

{\bf Proof of Theorem \ref{curve}.} By Convention \ref{conv3.1} we consider the case of curves in $\C^3$. Applying a polynomial automorphism to $X_2$ one can suppose that $S:=(X_1)_{\rm sing}$ coincides with
$(X_2)_{\rm sing}$ as a subset of $\C^n$.  Furthermore, by Proposition \ref{prop1.3} applying polynomial automorphisms to $X_i$ we 
can suppose that the restriction of any forgetting coordinate projection $p : \C^6_{x_1,y_1,z_1,x_2,y_2,z_2}\to \C^3$ to the graph $\Gamma$ of $\Psi$ yields
an isomorphism between $\Gamma$ and its image. By virtue of Remark \ref{rem1.1} (2) one can suppose that
these automorphisms do not ruin the equality $(X_1)_{\rm sing}=(X_2)_{\rm sing}$
since they can be chosen identical on $S$. Furthermore, applying a linear automorphism to both $X_1$ and $X_2$ we get the following:

(i) the restriction of any  coordinate to $S$ is injective;

(ii) the restriction of any coordinate to $X_i, i=1,2$ is a proper map onto $\C$,

(iii) for any any forgetting coordinate projection  $\C^6 \to \C^2$  its restriction to each $X_i$ yields a pseudo-isomorphism. 

Consider the the image $X_3$ (resp. $X_4$) of $\Gamma$ under the projection $\C^6 \to \C^3$ given by
$(x_1,y_1,z_1,x_2,y_2,z_2) \to (x_2,y_1,z_1)$ (resp. $(x_1,y_1,z_1,x_2,y_2,z_2) \to (x_2,y_2,z_1)$).
Suppose that $X_3$ and $X_4$ are naturally embedded in the same sample of $\C^3$ as $X_1$ and $X_2$.
By Proposition \ref{prop3.2} the natural isomorphism $X_1 \to X_3$ over $\C^2_{y_1,z_1}$ (resp.
$X_3 \to X_4$ over $\C^2_{x_2,z_1}$; resp. $X_4 \to X_2$ over $\C^2_{x_2,y_2}$) extends to a holomorphic
automorphism of $\C^3$. Taking composition of these holomorphic automorphisms we get the desired conclusion.
$\square$

\brem\label{rem3.2} We used only two specific properties of curves in the proof of Theorem \ref{curve}.
Namely, that 

(1)  the set  $S$ of singularities of $X_i\subset \C^n$, whose Zariski tangent space is of dimension $n$,  is finite and the dimension of
$T(X \setminus S)$ is at most $n-1$;

(2) for every
$x \in S$ each continuous meromorphic function vanishing at $x$ with sufficiently large multiplicity (independent on the function) is in fact holomorphic at $x$. 

Hence the same argument leads to a more general statement. 

{\em Let $\Psi : X_1 \to X_2$ be 
an isomorphism of closed affine subvarieties of $\C^{n}$ such that the singularities of $X_1$ satisfy conditions (1) and (2) above.
Suppose also that $n\geq 2 \dim X +1$.
Then $\Psi$ extends to a holomorphic automorphism of $\C^{n}$.}

\erem

\end{document}